\newif\ifpreprint
\preprinttrue  

\documentclass{article}

\usepackage{graphicx}
\usepackage{pgfplots}
\pgfplotsset{compat=1.18}
\usepackage{bm}
\usepackage{tikz}
\usepackage{amssymb,enumitem,mathtools}

\newenvironment{optproblem}[1][c]%
  {\setlength{\arraycolsep}{5pt}\begin{array}[#1]{ll}}%
  {\end{array}}
\usepackage{multirow}
\usepackage{tabularray}
\usepackage{csvsimple}
\usepackage{booktabs}
\usepackage{adjustbox}

\newcommand{\bnote}[1]{}
\newcommand{\snote}[1]{}
\newcommand{\anote}[1]{}
\newcommand{\nnote}[1]{}
\newcommand{\cnote}[1]{}

\newcommand{\eg}{{\it e.g.}}
\newcommand{\ie}{{\it i.e.}}

\newcommand{\ones}{\mathbf 1}
\newcommand{\reals}{{\mbox{\bf R}}}

\newcommand{\suppf}{{\phi}}

\newcommand{\tpose}{T}

\newcommand{\floor}[1]{\lfloor #1 \rfloor}
\newcommand{\ceil}[1]{\lceil #1 \rceil}

\newcommand{\myack}{Bartolomeo Stellato is supported by the NSF CAREER Award ECCS-2239771 and the ONR YIP Award N000142512147.
Stefano Gualandi acknowledges the contribution of the National Recovery and Resilience Plan, Mission 4 Component 2-Investment 1.4-National Center for HPC, Big Data and Quantum Computing (project code: CN\_00000013), funded by the European Union-NextGenerationEU.
The authors are thankful for NVIDIA for providing credits to access {\tt brev} cloud infrastructure with DGX B200 GPUs.}

\ifpreprint

\usepackage{fullpage}
\usepackage[colorlinks=true, citecolor=cyan, linkcolor=cyan, urlcolor=cyan]{hyperref}

\usepackage{natbib}

\author{Nicolas Blin$^{1}$, Stefano Gualandi$^{2}$, Christopher Maes$^{1}$, Andrea Lodi$^{3}$, Bartolomeo Stellato$^{4}$}

\date{\emph{\small NVIDIA$^1$} \\
\emph{\small Department of Mathematics, University of Pavia$^2$} \\
\emph{\small Jacobs Technion-Cornell Institute, Cornell Tech$^3$} \\
\emph{\small Department of Operations Research and Financial Engineering, Princeton University$^4$}\\[1em] 
\today}

\usepackage{amsmath, amsthm}

\theoremstyle{definition}

\theoremstyle{remark}

\newcommand{\keywords}[1]{\par\addvspace\baselineskip\noindent\textbf{Keywords:} #1}

\begin{document}

\title{Batched First-Order Methods for Parallel LP Solving in MIP}
\maketitle
\renewcommand{\qedsymbol}{$\blacksquare$}


\else



\usepackage{microtype}
\usepackage{graphicx}
\usepackage{subcaption}
\usepackage{booktabs} 

\usepackage{hyperref}

\newcommand{\theHalgorithm}{\arabic{algorithm}}

\usepackage{icml2026}



\usepackage{amsmath}
\usepackage{amssymb}
\usepackage{mathtools}
\usepackage{amsthm}

\usepackage[capitalize,noabbrev]{cleveref}

\usepackage{balance}

\theoremstyle{plain}
\newtheorem{theorem}{Theorem}[section]
\newtheorem{proposition}[theorem]{Proposition}
\newtheorem{lemma}[theorem]{Lemma}
\newtheorem{corollary}[theorem]{Corollary}
\theoremstyle{definition}
\newtheorem{definition}[theorem]{Definition}
\newtheorem{assumption}[theorem]{Assumption}
\theoremstyle{remark}
\newtheorem{remark}[theorem]{Remark}

\usepackage[textsize=tiny]{todonotes}
\renewcommand{\qedsymbol}{$\blacksquare$}

\begin{document}

\twocolumn[
  \icmltitle{Batched First-Order Methods for Parallel LP Solving in MIP}



  \icmlsetsymbol{equal}{*}

  \begin{icmlauthorlist}
    \icmlauthor{Firstname1 Lastname1}{equal,yyy}
    \icmlauthor{Firstname2 Lastname2}{equal,yyy,comp}
    \icmlauthor{Firstname3 Lastname3}{comp}
    \icmlauthor{Firstname4 Lastname4}{sch}
    \icmlauthor{Firstname5 Lastname5}{yyy}
    \icmlauthor{Firstname6 Lastname6}{sch,yyy,comp}
    \icmlauthor{Firstname7 Lastname7}{comp}
    \icmlauthor{Firstname8 Lastname8}{sch}
    \icmlauthor{Firstname8 Lastname8}{yyy,comp}
  \end{icmlauthorlist}

  \icmlaffiliation{yyy}{Department of XXX, University of YYY, Location, Country}
  \icmlaffiliation{comp}{Company Name, Location, Country}
  \icmlaffiliation{sch}{School of ZZZ, Institute of WWW, Location, Country}

  \icmlcorrespondingauthor{Firstname1 Lastname1}{first1.last1@xxx.edu}
  \icmlcorrespondingauthor{Firstname2 Lastname2}{first2.last2@www.uk}

  \icmlkeywords{Machine Learning, ICML}

  \vskip 0.3in
]

\icmlkeywords{Mixed-Integer Linear Programming, First-order method, Parallel computing, GPU architectures}

\printAffiliationsAndNotice{\icmlEqualContribution}

\fi

\begin{abstract}
  We present a batched first-order method for solving multiple linear programs in parallel on GPUs. Our approach extends the primal-dual hybrid gradient algorithm to efficiently solve batches of related linear programming problems that arise in mixed-integer programming techniques such as strong branching and bound tightening. By leveraging matrix-matrix operations instead of repeated matrix-vector operations, we obtain significant computational advantages on GPU architectures. 
  We demonstrate the effectiveness of our approach on various case studies and identify the problem sizes where first-order methods outperform traditional simplex-based solvers depending on the computational environment one can use.
  This is a significant step for the design and development of integer programming algorithms tightly exploiting GPU capabilities where we argue that some specific operations should be allocated to GPUs and performed in full instead of using light-weight heuristic approaches on CPUs.
\end{abstract}

\ifpreprint
\fi

\section{Introduction}
\label{sec:intro}

First-order methods for solving linear programming (LP) problems have recently been rediscovered for the (approximate) solution of very large instances.
This is due to the fact that first-order methods are a natural choice for parallelization, which in turn is the relevant characteristic to look into for exploiting GPUs.

Although there has been significant progress in solving LPs with first-order methods starting from the first GPU implementations of operator-splitting solvers based on the alternating direction method of multipliers (ADMM), such as the SCS solver~\cite{ocpb:16,odonoghue:21} and the OSQP solver~\cite{stellatoOSQPOperatorSplitting2020,schubiger2020gpu} and, more recently, the primal-dual hybrid gradient (PDHG) algorithm~\cite{chambolle2011first} implementation in the PDLP solver~\cite{Applegate2021,Applegate2023,lu2025cupdlp}, the target of this research area has been (so far) solving extremely large LPs and there is some skepticism that this technology can be applied in mixed-integer programming (MIP)~\cite{Rothberg2024}. This is due to the fact that interior point algorithms are also very effective and parallelizable (so, solving the LP relaxation of any MIP is somehow doable) and the simplex method seems unbeatable in its ability of warm-starting in the presence of a small modification of the solution space due to the change of a variable bound or the addition of a round of cuts.

The GPU use in MIP has been also limited by the technological burden of moving large amounts of data from CPUs (where the computation has been run traditionally) to GPUs and vice versa. This difficulty has been observed for example in the context of the use of machine learning (ML) models to augment branch-and-bound methods \cite{scavuzzo2024machine}. Indeed, some significant success has been already obtained (see, \eg, \cite{bonami2022classifier,berthold2025cuts}), but the CPU vs GPU interaction has slowed down the adoption of ML models using neural networks (NN) representations, especially when those models need to be run at each node of the branch-and-bound tree like, for example, for variable selection in branching \cite{Gasse2019,gupta2020hybrid}.

In other words, it is possible that an effective GPUs exploitation in MIP -- besides the  technological improvement in GPU memory -- requires rethinking the branch and bound itself, potentially abandoning some of its foundational certainties (the simplex method) and removing its computational work limits. 

Along this rethinking path, the current paper shows that one GPU benefit for MIP is on exploiting the parallelization power to execute in full on GPUs chunks of the computation that, traditionally, are performed by clever heuristics that guarantee good approximations and limit the computational load. This is the case of strong branching and optimization-based bound tightening that involve the solution of batches of LPs with minimal differences (variable bounds). 
State-of-the-art implementations of batched LP solutions arise in differentiable optimization, where problems are structured as neural network layers \cite{cvxpylayers2019,lu2024mpax,besancon2023diffopt}, and primarily use basic vectorization operations (\eg, {\tt vmap}). These methods do not fully exploit the parallelism of GPU-based matrix operations and require data duplication across batch instances since many solvers lack full thread safety for parallel instance solving.

In this paper, we propose a fully parallelizable first-order method directly designed to solve batches of LPs all at once on GPUs, overcoming limitations of existing batch approaches. In the case of strong branching, this allows a precise initialization of pseudocosts (the capillary information MIP solvers use for guiding branching) in only a few rounds of memory communication, thus gaining quality without paying a high price for GPU vs CPU interaction. 

The remainder of the paper is organized as follows. In Section \ref{sec:literature}, we review the relevant literature concerning first-order methods for LPs and their GPU implementations. In Section \ref{sec:formulation}, we present the way we formulate an LP so as to be amenable to the first-order algorithm described in Section \ref{sec:PDHG}. In Section \ref{sec:batchLPs}, we discuss our contribution in solving batches of LPs by first-order methods. In Section \ref{sec:results}, we present our computational experiments with our implementation of {\sc BatchLP} to perform strong branching and OBBT on standard benchmarks from the literature. Finally, in Section \ref{sec:conclusions}, we draw some conclusions.

%

\section{Literature review}
\label{sec:literature}

The first effective first-order method for LP is the PDLP algorithm \cite{Applegate2021,Applegate2023}. PDLP applies the primal-dual hybrid gradient (PDHG) method, see \cite{chambolle2011first}, and uses restarting and averaging to accelerate it. Several other enhancements were included in PDLP such as adaptive step-size (see also \cite{Chambolle2024}) and problem scaling.

Building on PDLP, the Halpern Peaceman-Rachford (HPR) method takes a weighted average between the current PDHG iterate and an initial point \cite{Lu2024}. The initial point is updated when certain restart conditions are met. This leads to an algorithm called restarted Halpern PDHG (rHPDHG). A similar restarted Halpern algorithm is developed in \cite{Chen2024} with a relaxation step that results in a longer step size. Inspired by \cite{Chen2024}, an extension is considered in \cite{Lu2024}, called reflected restarted Halpern PDHG (r$^2$HPDHG), where the Halpern iteration is performed on the reflection of the PDHG operator instead of the operator itself. Finally, the extension of rHPDHG to conic LP is presented in \cite{Xiong2024}.

In terms of successive improvements of the PDLP algorithm, a geometric interpretation of the PDLP behavior is presented in \cite{Liu2024} and based on such interpretation a new crossover algorithm is designed to recover a vertex solution for an LP. A key role in this type of algorithms is played by infeasibility detection, \ie, recognizing infeasible subproblems, so as to avoid useless iterations. This is extensively explored in \cite{Applegate2024,Banjac2019}. 

The idea of extending a first-order method to solve batches of linear programs on the GPU for strong branching was first introduced in \cite{Nair2020}. The goal and methodology in \cite{Nair2020} are fundamentally different from our proposal because the GPU implementation of strong branching is used to collect data for training an ML model at the branching task, \ie, to approximate strong branching. Instead, in this paper, we propose to use GPUs to do strong branching, not approximate it, for a certain number of iterations that allow to initialize pseudocosts. Besides the different motivation, the way in which the GPUs are used for strong branching is very different. In \cite{Nair2020}, the batches of LPs are solved by a modification of ADMM, and in a setting where instead of branching on fractional variables the MIP algorithm branches on all variables. Finally, the assessment of the speed-up with respect to traditional CPU implementations does not consider the fact that each LP can be solved by starting from the optimal basis of the simplex method, so generally doing only very few pivots to converge. This leads to an overestimation of the speed-up while our methodology assesses it very carefully in different benchmark settings, see Section \ref{sec:results}. 

\section{Problem formulation}
\label{sec:formulation}

Consider a linear optimization problem in the primal-dual form
\begin{equation}\label{eq:main_lp}
\ifpreprint\else\!\!\!\!\!\!\!\fi\begin{optproblem}
    \text{min} & c^\tpose x\\
    \text{s.t.} & l \le Ax \le u,\\ &\underline{x} \le x \le \bar{x},
\end{optproblem}
\; 
\begin{optproblem}
    \text{max} &  -\suppf_{[\underline{x}, \bar{x}]}(r) -\suppf_{[l, u]}(y)\\
    \text{s.t.} & c + A^\tpose y + r = 0\\
    &y \in B_{[l, u]},\quad r \in B_{[\underline{x}, \bar{x}]},
\end{optproblem}
\end{equation}
where the decision variables are $x \in \reals^n$, and the coefficients of the linear objective function are $c \in \reals^n$.
The constraints are defined by matrix $A \in \reals^{m \times n}$ and vectors $l$, $u$, $\underline{x}$, and $\bar{x}$, with $l_i ,\underline{x}_i\in \{- \infty\} \cup \reals$ and $u_i,\bar{x}_i \in \reals \cup \{+\infty\}$. We represent equality constraints as $l_i = u_i$ and variable fixings as $\underline{x}_i=\bar{x}_i$.
We define the dual variables for the inequality constraints as~$y \in \reals^m$ and for variable bounds as~$r \in \reals^n$.
The set $B$ represents the barrier cone of a hyperrectangle, \ie, $B_{[a, b]} = \{v \mid v_i \ge 0\:\text{if}\: a_i=-\infty;\; v_i \le 0\:\text{if}\: b_i=\infty;\;v_i \in \reals \;\text{otherwise}\}$.
Function $\phi$ represents the support function of a hyperrectangle, \ie, $\phi_{[a, b]}(v) = \sup_{z \in [a, b]} v^\tpose z = b^\tpose v_{+} + a^\tpose v_{-}$, where 
$v_{+} = \max\{v, 0\}$ and $v_{-} = \min\{v, 0\}$.


To derive our algorithm, we reformulate~\eqref{eq:main_lp} as the following saddle-point problem
\begin{equation}\label{eq:lp_saddle}
    \max_{y} \min_{\underline{x} \le x \le \bar{x}} c^\tpose x + y^\tpose Ax - \suppf_{[l, u]}(y).
\end{equation}


\section{Primal-dual hybrid gradient to solve a single instance}
\label{sec:PDHG}

We apply the primal-dual hybrid gradient method (PDHG)~\cite{chambolle2011first} with reflected Halpern iterations~\cite{halpern1967fixed,Lu2024,cupdlpx}
to solve problem~\eqref{eq:main_lp}. 
The iterations in terms of the primal-dual pair $z^k = (x^k, y^k) \in \reals^{n + m}$ consist of
\begin{equation}\label{eq:mainiter}
    \begin{array}{ll}
         z^{k+1} &= 
\dfrac{k+1}{k+2} \left( 2 T(z^k) 
- z^k \right) 
+ \dfrac{1}{k+2} z^0,
    \end{array}
\end{equation}
with the main operator $T$ being defined as
\ifpreprint
\begin{equation*}
     T(z^k) = \left\{ (x, y) \;\middle|\;
          \begin{array}{@{}l@{}}
          x = \Pi_{[\underline{x}, \bar{x}]}\left( x^k - \tau (c + A^\tpose y^k)\right)\\[4pt]
          y = y^k + \sigma A(2x - x^{k}) - \sigma \Pi_{[l, u]} ( \sigma^{-1} y^k + A(2x - x^{k}) )
         \end{array}
    \right\}.
\end{equation*}
\else
\begin{equation*}
     T(z^k) = \left\{ (x, y) \;\middle|\;
          \begin{array}{@{}l@{}}
          x = \Pi_{[\underline{x}, \bar{x}]}\left( x^k - \tau (c + A^\tpose y^k)\right)\\[4pt]
          y = y^k + \sigma A(2x - x^{k}) \\
          \quad {-}\, \sigma \Pi_{[l, u]} ( \sigma^{-1} y^k + A(2x - x^{k}) )
         \end{array}
    \right\}.
\end{equation*}
\fi
Here, $z^0$ is the initial iterate used in the Halpern iterations as the anchor point. Moreover, $\Pi_{[a, b]}$ is the Euclidean projection onto a hypercube defined as the elementwise operation $\Pi_{[a, b]}(v) = \max\{\min\{v, b\}, a\}$.
The primal and dual step-sizes are~$\tau$ and~$\sigma$, respectively.
As commonly done in PDHG~\cite{Applegate2021}, we parametrize the step sizes as $\tau = \eta/w$ and $\sigma = \eta w$, where we refer to $\eta$ as the step size and $w$ as the primal weight.

\paragraph{Restarts.}
Similarly to~\cite{cupdlpx}, we apply an adaptive restart scheme based on the fixed-point residual progress of the non-reflected iterates. 
More specifically, we define the fixed-point residual metric 
\begin{equation}\label{eq:residual}
    r(z) = \|T(z) - z\|_M,\, \text{with}\,  M = \begin{bmatrix} (w/\eta) I & A^\tpose\\ A & (1/(\eta w))I\end{bmatrix}.
\end{equation}
To restart, we partition the iterations in two loops, an outer loop indexed by $n$ and an inner loop indexed by $k$, with corresponding iterates $z^{n, k}$.
The algorithm repeats iteration~\eqref{eq:mainiter} until one of the following restart conditions is met: 
\begin{itemize}[noitemsep]
    \item sufficient decay: $r(z^{n, k}) \le \beta_s r(z^{n, 0})$,
    \item necessary decay and no inner progress: $r(z^{n, k}) \le \beta_n r(z^{n, 0})$ and $r(z^{n, k}) > r(z^{n, k-1})$, and
    \item iteration limit: $k > \beta_a K$, with $K$ being the total number of iterations.
\end{itemize}
In this case, the anchor point $z^0$ is set to the current iterate.

\paragraph{Step-size and primal weight.}
We adopt the constant step-size from~\cite{cupdlpx}, where $\eta = 0.998/\|A\|_2$ and the exponential smoothing technique to update the primal weight $w$~\cite{applegate2025pdlppracticalfirstordermethod}.
Specifically, whenever $d^n = \|x^{n+1, 0} - x^{n, 0}\|/\|y^{n+1, 0} - y^{n, 0}\|$ is finite at every outer loop iteration, we update $w^{n+1} = \exp\left(\theta \log(d^n) + (1 - \theta) \log(w^n)\right)$.

\paragraph{Stopping criteria.}
This algorithm, including restarts, has been analyzed in~\cite{Lu2024} and shown to converge at a linear rate. We terminate the algorithm if the following optimality conditions are satisfied:
\ifpreprint
\begin{equation}\label{eq:stopping}
    \begin{array}{@{}l@{}}
         |c^\tpose x^k + \suppf_{[\underline{x}, \bar{x}]}(r^k) + \suppf_{[l, u]}(y^k)| \le \epsilon (1 + |c^\tpose x^k| + |\suppf_{[\underline{x}, \bar{x}]}(r^k) + \suppf_{[l, u]}(y^k)|)\\[4pt]
         \|Ax^k - \Pi_{[l, u]}(Ax^k) \|_2 \le \epsilon (1 + \|Ax^k\|_2)\\[4pt]
         \|c + A^\tpose y^k + r^k\|_2 \le \epsilon (1 + \|c\|_2).
    \end{array}
\end{equation}
\else
\begin{equation}\label{eq:stopping}
    \begin{array}{@{}l@{}}
         |c^\tpose x^k + \suppf_{[\underline{x}, \bar{x}]}(r^k) + \suppf_{[l, u]}(y^k)|\\
         \quad \quad \le \epsilon (1 + |c^\tpose x^k| + |\suppf_{[\underline{x}, \bar{x}]}(r^k) + \suppf_{[l, u]}(y^k)|)\\[4pt]
         \|Ax^k - \Pi_{[l, u]}(Ax^k) \|_2 \le \epsilon (1 + \|Ax^k\|_2)\\[4pt]
         \|c + A^\tpose y^k + r^k\|_2 \le \epsilon (1 + \|c\|_2).
    \end{array}
\end{equation}
\fi
Since iteration~\eqref{eq:mainiter} does not track a reduced cost vector $r^k \in \reals^n$, we compute it as~$r^k = \Pi_{B_{[\underline{x}, \bar{x}]}}(-c - A^\tpose y^k)$ where $\Pi_{B_{[\underline{x}, \bar{x}]}}$ is the Euclidean projection on the barrier cone of the variable bounds.
In the actual implementation, we replace $\Pi_{B_{[\underline{x}, \bar{x}]}}$ with an alternative formulation that tends to be more robust to bounds taking large values, see Appendix A of~\cite{applegate2025pdlppracticalfirstordermethod}. 

\paragraph{Infeasibility detection.}
 To detect infeasibility, we apply the  conditions recently developed based on the infimal displacement vector of operator splitting algorithms~\cite{pazy1971asymptotic,Banjac2019,Applegate2024}, specifically applied to Halpern iterations~\cite{pmlr-v202-park23k}.
 In particular, given the primal-dual iterates~$x^{k-1}, y^{k-1}$ and reduced cost~$r^{k-1}$, we apply the operator~$T$ obtaining $(\tilde{x}^{k}, \tilde{y}^{k}) = T(x^{k-1}, y^{k-1})$ and $\tilde{r}^{k} = \Pi_{B_{[\underline{x}, \bar{x}]}}(-c - A^\tpose\tilde{y}^{k})$. 
 Then, we compute the displacement vectors $\delta x^k = \tilde{x}^{k} - x^{k-1}$, $\delta y^k = \Pi_{B_{[l, u]}}(\tilde{y}^{k} - y^{k-1})$, and $\delta r^k = \Pi_{B_{[\underline{x}, \bar{x}]}}(\tilde{r}^{k} - r^{k-1})$. We define the $\epsilon$-approximate primal infeasibility conditions as
 \begin{equation}\label{eq:primal-infeas}
 \ifpreprint\else\!\!\!\!\!\!\fi\begin{array}{ll}
     &\phi_{[l, u]}(\delta y^k) + \phi_{[\underline{x}, \bar{x}]}(\delta r^k) < 0, \\ 
     & \left\|A^\tpose \delta y^k + \delta r^k\right\| \le \epsilon \left(\phi_{[l, u]}(\delta y^k) + \phi_{[\underline{x}, \bar{x}]}(\delta r^k)\right).
     \end{array}
 \end{equation}
Similarly, the $\epsilon$-approximate dual infeasibility conditions are
 \begin{equation}\label{eq:dual-infeas}
 \begin{array}{ll}
     c^\tpose \delta x^k < 0, \; \left\|\delta x^k - \Pi_{R_{[\underline{x}, \bar{x}]}}(\delta x^k)\right\| \le \epsilon |c^\tpose \delta x^k|, \\ \left\|A\delta x^k - \Pi_{R_{[l, u]}}(A\delta x^k)\right\| \le \epsilon |c^\tpose \delta x^k|,\\    
     \end{array}
 \end{equation}
where $\Pi_{R_{[a, b]}}$ is the Euclidean projection onto the recession cone of hyperrectangle $[a, b]$, which is defined as 
$R_{[a, b]} = \{v \mid v_i \ge 0\:\text{if}\: b_i=\infty;\; v_i \le 0\:\text{if}\: l_i=-\infty;\;v_i =0 \;\text{otherwise}\}$.

\section{Solving batches of problems}
\label{sec:batchLPs}


Instead of solving a single instance of problem~\eqref{eq:main_lp}, we consider solving a \emph{batch} of $N$ instances with the same data matrix $A$, and with varying problem vectors. 
We stack the problem vectors in the objective matrix $C = [\,c^1\;\cdots\; c^N\,] \in \reals^{n\times N}$ and the bounds matrices $\underline{X} = [\,\underline{x}^1\;\cdots\; \underline{x}^N\,]$, $\overline{X}=[\,\bar{x}^1\;\cdots\; \bar{x}^N\,] \in \reals^{n \times N}$, $L = [\,l^1\;\cdots\; l^N\,]$, and $U=[\,u^1\;\cdots\; u^N\,] \in \reals^{n \times N}$.
By defining the matrices of primal iterates as $X^k = [\,x_1^k\;\cdots\;x_N^k\,]\in \reals^{n\times N}$, dual iterates as~$Y^k = [\,y_1^k\;\cdots\;y_N^k\,]\in \reals^{m\times N}$, and combined iterates as~$Z^k = [\,z_1^k\;\cdots\;z_N^k\,]\in \reals^{(n+m)\times N}$ we can rewrite iteration~\eqref{eq:mainiter} as
\begin{equation}\label{eq:batched_iter}
    \begin{array}{ll}
         Z^{k+1} &= 
\dfrac{k+1}{k+2} \left( 2 T(Z^k) 
- Z^k \right) 
+ \dfrac{1}{k+2} Z^0,
    \end{array}
\end{equation}
with the operator $T(Z^k) = (X, Y)$ defined componentwise as
\ifpreprint
\begin{equation}\label{eq:main-iterate}
          \left\{\begin{array}{@{}l@{}}
          X = \Pi_{[\underline{X}, \overline{X}]}\left( X^k - \tau \otimes (C + A^\tpose Y^k)\right)\\[4pt]
          Y = Y^k + \sigma \otimes A(2X - X^{k}) - \sigma \otimes \Pi_{[L,U]} ( \sigma^{-1} \otimes Y^k + A(2X - X^{k}) ),
         \end{array}\right.
\end{equation}
\else
\begin{equation}\label{eq:main-iterate}
          \left\{\begin{array}{@{}l@{}}
          X = \Pi_{[\underline{X}, \overline{X}]}\left( X^k - \tau \otimes (C + A^\tpose Y^k)\right)\\[4pt]
          Y = Y^k + \sigma \otimes A(2X - X^{k}) \\
          \quad \quad {-}\, \sigma \otimes \Pi_{[L,U]} ( \sigma^{-1} \otimes Y^k + A(2X - X^{k}) ),
         \end{array}\right.
\end{equation}
\fi
\noindent and the projection operators representing the elementwise matrix projection.
In addition, we have $w, \tau, \sigma \in \reals^N$, since we keep a different value for each subproblem.
As for the notation of primal weight and step sizes, we define the Kronecker products~$\sigma \otimes X := [\sigma_1 x_1, \dots, \sigma_N x_N]$ and $\tau \otimes Y := [\tau_1 y_1, \dots, \tau_N y_N]$, that is, every column of $X$ and $Y$ gets a different step size.

\paragraph{Restarts.}
To adapt the restarting logic to take care of batches of LPs, we define the \emph{averaged} fixed-point residual metric as $\tilde{r}(Z^k) = (1/N)\sum_{j=1}^N r(z_j^k)$, 
that is, the average residual of the columns of matrix $Z^k$, where each residual is computed as in~\eqref{eq:residual}, but using the corresponding values of primal weight~$w_j$ and step-sizes~$\tau_j$ and~$\sigma_j$.
Whenever the average residual $\tilde{r}(Z^k)$ satisfies one of the three conditions as the ones for the single-instance residual~\eqref{eq:residual},
we restart all the batches at the same time, by recomputing the new primal weights~$w_j^{n+1}$ by exponential smoothing, and the step-sizes~$\tau_j^{n+1}$ and $\sigma_j^{n+1}$ accordingly, for each problem~$j=1,\dots, N$. 

\paragraph{Stopping criteria and infeasibility detection.}
The stopping conditions for our batched algorithm are as in \eqref{eq:stopping}, but rewritten in matrix notation.
The same apply to infeasibility criteria \eqref{eq:primal-infeas} and \eqref{eq:dual-infeas}. Whenever we finish solving a problem in the batch, we permute the iterate matrices so that the final columns correspond to the solved problems; equivalently, we keep track of $\pi = (\pi_1, \dots, \pi_N)$ as a permutation of the indices, with $X^k = [\,x_{\pi_1}^k\;\cdots\;x_{\pi_N}^k\,]$ and~$Y^k = [\,y_{\pi_1}^k\;\cdots\;y_{\pi_N}^k\,]$.
 In this way, we stop updating those variables and reduce the kernel sizes dynamically.

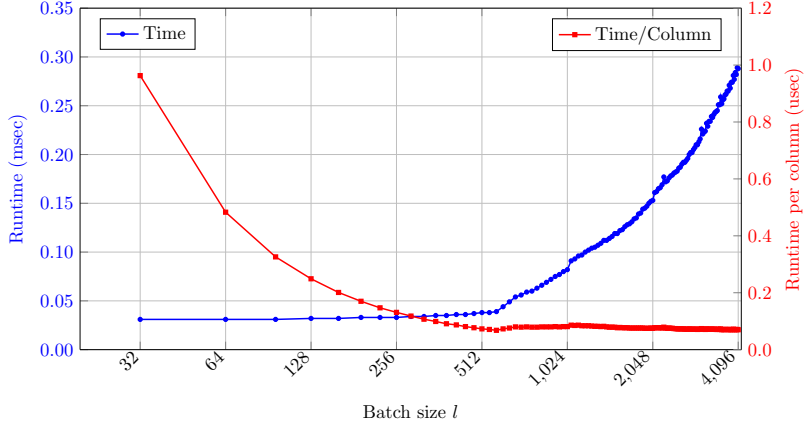
\begin{figure*}[t]
    \centering
    \begin{adjustbox}{width=\ifpreprint 0.65\textwidth\else 0.8\textwidth\fi}

\begin{tikzpicture}
\begin{axis}[
    xmode=log,
    log basis x=2,
    xticklabel={\pgfmathparse{2^(\tick)}\pgfmathprintnumber[fixed,precision=0]{\pgfmathresult}},
    width=14cm,
    height=8cm,
    xlabel={Batch size $l$},
    ymin=0, ymax=0.35,
    ylabel style={blue},
    ylabel={Runtime (msec)},
    yticklabel style={blue},
    yticklabel style={/pgf/number format/.cd, fixed, precision=2, zerofill},
    axis y line*=left,
    xmin=0, xmax=4200,
    xtick={32,64,128,256,512,1024,2048,4096},
    xticklabel style={rotate=45, anchor=east},
    grid=major,
    legend pos=north west,
]
\addplot[blue, thick, mark=*, mark size=1pt] table[col sep=comma] {
Size,Total
32,0.031
64,0.031
96,0.031
128,0.032
160,0.032
192,0.033
224,0.033
256,0.033
288,0.034
320,0.034
352,0.035
384,0.035
416,0.036
448,0.036
480,0.037
512,0.038
544,0.038
576,0.039
608,0.044
640,0.049
672,0.054
704,0.056
736,0.059
768,0.060
800,0.063
832,0.066
864,0.069
896,0.072
928,0.075
960,0.077
992,0.080
1024,0.082
1056,0.091
1088,0.093
1120,0.096
1152,0.097
1184,0.100
1216,0.102
1248,0.104
1280,0.105
1312,0.107
1344,0.109
1376,0.112
1408,0.112
1440,0.114
1472,0.116
1504,0.119
1536,0.119
1568,0.122
1600,0.123
1632,0.126
1664,0.128
1696,0.129
1728,0.131
1760,0.134
1792,0.135
1824,0.139
1856,0.140
1888,0.144
1920,0.145
1952,0.147
1984,0.150
2016,0.152
2048,0.153
2080,0.161
2112,0.162
2144,0.165
2176,0.166
2208,0.169
2240,0.177
2272,0.172
2304,0.173
2336,0.176
2368,0.178
2400,0.179
2432,0.181
2464,0.182
2496,0.183
2528,0.186
2560,0.187
2592,0.190
2624,0.192
2656,0.192
2688,0.194
2720,0.196
2752,0.200
2784,0.202
2816,0.202
2848,0.205
2880,0.207
2912,0.210
2944,0.210
2976,0.213
3008,0.216
3040,0.226
3072,0.221
3104,0.223
3136,0.224
3168,0.232
3200,0.229
3232,0.234
3264,0.234
3296,0.239
3328,0.238
3360,0.241
3392,0.243
3424,0.244
3456,0.245
3488,0.251
3520,0.251
3552,0.259
3584,0.252
3616,0.256
3648,0.257
3680,0.261
3712,0.262
3744,0.265
3776,0.265
3808,0.271
3840,0.268
3872,0.274
3904,0.274
3936,0.281
3968,0.277
4000,0.284
4032,0.282
4064,0.289
4096,0.288
};
\legend{Time}

\end{axis}

\begin{axis}[
    xmode=log,
    log basis x=2,
    xticklabel={\pgfmathparse{2^(\tick)}\pgfmathprintnumber[fixed,precision=0]{\pgfmathresult}},
    width=14cm,
    height=8cm,
    ylabel={Runtime per column (usec)},
    ylabel style={red},
    yticklabel style={red},
    yticklabel style={/pgf/number format/.cd, fixed, precision=1, zerofill},
    axis y line*=right,
    axis x line=none,
    xmin=0, xmax=4200,
    ymin=0, ymax=1.2,
    legend pos=north east,
]

\addplot[red, thick, mark=square*, mark size=1pt] table[col sep=comma] {
Size,TotalCol
32,0.963
64,0.483
96,0.326
128,0.249
160,0.201
192,0.170
224,0.147
256,0.131
288,0.118
320,0.107
352,0.099
384,0.091
416,0.087
448,0.081
480,0.077
512,0.073
544,0.071
576,0.068
608,0.073
640,0.076
672,0.080
704,0.079
736,0.080
768,0.079
800,0.079
832,0.080
864,0.080
896,0.080
928,0.081
960,0.080
992,0.081
1024,0.081
1056,0.086
1088,0.085
1120,0.086
1152,0.084
1184,0.084
1216,0.084
1248,0.083
1280,0.082
1312,0.082
1344,0.081
1376,0.082
1408,0.080
1440,0.079
1472,0.079
1504,0.079
1536,0.077
1568,0.078
1600,0.077
1632,0.077
1664,0.077
1696,0.076
1728,0.076
1760,0.076
1792,0.076
1824,0.076
1856,0.076
1888,0.076
1920,0.075
1952,0.075
1984,0.075
2016,0.076
2048,0.075
2080,0.077
2112,0.077
2144,0.077
2176,0.077
2208,0.076
2240,0.079
2272,0.076
2304,0.075
2336,0.076
2368,0.075
2400,0.075
2432,0.074
2464,0.074
2496,0.073
2528,0.073
2560,0.073
2592,0.073
2624,0.073
2656,0.072
2688,0.072
2720,0.072
2752,0.073
2784,0.072
2816,0.072
2848,0.072
2880,0.072
2912,0.072
2944,0.071
2976,0.072
3008,0.072
3040,0.074
3072,0.072
3104,0.072
3136,0.071
3168,0.073
3200,0.072
3232,0.072
3264,0.072
3296,0.072
3328,0.072
3360,0.072
3392,0.072
3424,0.071
3456,0.071
3488,0.072
3520,0.071
3552,0.073
3584,0.070
3616,0.071
3648,0.070
3680,0.071
3712,0.071
3744,0.071
3776,0.070
3808,0.071
3840,0.070
3872,0.071
3904,0.070
3936,0.072
3968,0.070
4000,0.071
4032,0.070
4064,0.071
4096,0.070
};
\legend{Time/Column}
\end{axis}
\end{tikzpicture}
\end{adjustbox}
\caption{Sparse Matrix-Matrix Multiplies vs Batch Size on \texttt{csched007}. Total time for 10 $AX$ and 10 $A^\tpose Y$ operations (blue circle markers on the left y-axis) and time/column (red square markers on the right y-axis).}    
\label{fig:csched007}
\end{figure*}

\paragraph{Matrix-matrix products and optimal batch size.}
The major bottleneck in the batched algorithm are matrix-matrix products $A^\tpose Y$ and $AX$, which can be significantly accelerated on GPUs.
Moreover, observe that the cost of computing $A X^k$ and $A^\tpose Y^k$ is weakly dependent on the batch size (or the number of columns in $X^k$ and $Y^k$). We varied the batch size and timed 10 matrix-matrix products of $A X^0$ and 10 matrix-matrix products of $A^\tpose Y^0$ using the constraint matrix $A$ from the MIPLIB benchmark problem \texttt{csched007} \cite{Yunes2009}. 
Figure~\ref{fig:csched007} shows the times that CUDA library cuSPARSE takes to compute these matrix-matrix products on a RTX 4500 ADA GPU. Observe that from batch size 32 to 512 the time taken for these products remains constant at approximately 0.03 ms. 
After a batch size of 512 the time for these products increases with the batch size. 
We define the optimal batch size as the one that maximizes the number of matrix-vector multiplies per second. 
The optimal batch size, influenced by factors like the constraint matrix $A$, GPU, and software configuration, cannot be predetermined. 
For a given problem and GPU, we estimate it by timing a few matrix-matrix products before solving the MIP.   

\paragraph{Naive batched implementation.}
Although the computational speed-ups achievable by matrix–matrix multiplies over batches seem very promising, a naive batched implementation of our algorithm could be very inefficient due to data movement between CPU and GPU, which causes significant memory overhead. In the next section we tailor our batched first-order method to strong branching, a core task in mixed integer programming. In Appendix~\ref{app:obbt}, we outline similar customizations to implement optimization-based bound tightening.

\subsection{Full strong branching}
\label{sec:strongbranching}

A key step in branch-and-bound algorithms for MIP is choosing which fractional variable to branch on at each node.
Full strong branching (FSB) solves two LPs for each fractional variable: one that rounds up and one that rounds down~\cite{Achterberg2005}.
The resulting objective values determine the branching score.

Precisely, consider an optimal solution $x^{\rm rel}$ of the LP relaxation at the current node with $p$ fractional variables.
For notational simplicity, assume the first $p$ components of $x^{\rm rel}$ are fractional.
FSB solves $N = 2p$ LPs: for each fractional variable $x_i$ ($i = 1, \ldots, p$), we solve one LP with $x_i \geq \ceil{x_i^{\rm rel}}$ (branch up) and one with $x_i \leq \floor{x_i^{\rm rel}}$ (branch down).
All $N$ problems share the same cost vector $c$, constraint bounds $l$, $u$, and base variable bounds $\underline{x}$, $\bar{x}$; only one bound component differs per problem.

In matrix notation, the objective and constraint bounds are $C = c\ones^\tpose$, $L = l\ones^\tpose$, and $U = u\ones^\tpose$.
The variable bound matrices $\underline{X}, \overline{X} \in \reals^{n \times N}$ have columns equal to $\underline{x}$ and $\bar{x}$, respectively, except
\begin{itemize}[noitemsep,nosep]
    \item Column $i$ (branch up on $x_i$): $(\underline{X})_{ii} = \ceil{x_i^{\rm rel}}$,
    \item Column $p + i$ (branch down on $x_i$): $(\overline{X})_{i,p+i} = \floor{x_i^{\rm rel}}$.
\end{itemize}

Since only one bound entry changes per problem, we avoid forming $\underline{X}$ and $\overline{X}$ explicitly; instead, we index into $x^{\rm rel}$, $\underline{x}$, and $\bar{x}$ as needed.
Once memory is allocated for $p$ fractional variables, the same allocation serves smaller values of $p$, which naturally occurs as the tree deepens and more variables become fixed.

\subsection{Optimization-based bound tightening}
\label{app:obbt}
Optimization-based bound tightening (OBBT) is a core preprocessing step in MIP solvers, in which we update the bounds of each variable as $\underline{x}_i = \max\{\underline{x}_i, \underline{x}_i^{\rm obbt}\}$ and $\bar{x}_i = \min\{\bar{x}_i, \bar{x}_i^{\rm obbt}\}$, where 
\begin{equation*}
    \underline{x}_i^{\rm obbt} = \begin{optproblem}[t]
        \text{min} & x_i\\
        \text{s.t.} & l \leq A x \leq u\\
        & \underline{x} \leq x \leq \bar{x}
    \end{optproblem}\;
    \bar{x}_i^{\rm obbt} = \begin{optproblem}[t]
        \text{max} & x_i\\
        \text{s.t.} & l \leq A x \leq u\\
        & \underline{x} \leq x \leq \bar{x},
    \end{optproblem}
\end{equation*}
for $i=1,\dots,n$.
This corresponds to solving $2n=N$ LPs having the same constraint bounds, $l, u$, and variable bounds $\underline{x},\bar{x}$.
The objectives correspond to the maximization and minimization of each component of $x$.
In matrix form, this corresponds to matrix bounds~$L=l\ones^\tpose, U=u\ones^\tpose$, $\underline{X} = \underline{x} \ones^\tpose$, and $\overline{X} = \overline{x} \ones^\tpose$, all in $\reals^{n \times N}$.
The objective matrix becomes~$C = [\,I\;-I\,] \in \reals^{n \times N}$.
Similarly to FSB, we never materialize matrix $C$, where we instead smartly reference its columns that correspond to the canonical basis vectors with alternating signs.\\

Before ending Section \ref{sec:batchLPs}, it is important to note that none of the MIP solvers, either commercial or academic, uses FSB or performs OBBT extensively on all variables. A number of work limits are used to maintain the computational footprint of those methods — which are widely acknowledged to be extremely powerful — under control. For strong branching, one selects a restricted candidate subset of variables and the LPs are solved by limiting the number of Simplex pivots. For OBBT, the LPs of some variables are solved (heuristically) once triggered by conditions capturing the likelihood of strengthening the variable bound. In other words, the computational section below aims at answering the question: can one exploit FSB and OBBT at their full strength by moving the associated computation to GPUs with our proposed algorithmic approach?

\section{Numerical experiments}
\label{sec:results}


In this section, we present our computational results for full strong branching and OBBT using benchmark instances from the literature. Specifically, for both strong branching and OBBT, we evaluate the performance of our batched algorithms described in Sections~\ref{sec:strongbranching} and \ref{app:obbt}, respectively, by solving LPs as if we were at the root node of a branch-and-bound tree.
\begin{table*}[t]
\centering
\caption{Runtime (in seconds) comparison of dual simplex and {\sc BatchLP} as implemented in cuOpt 26.02.\label{tab:cuopt:threads}}
\ifpreprint\footnotesize\else\small\fi
\setlength{\tabcolsep}{3.5pt}
\begin{tabular}{lrrrrrrrrr}
\toprule
 & \multicolumn{4}{c}{Problem} & \multicolumn{4}{c}{cuOpt 26.02 -- Dual Simplex} & \multicolumn{1}{c}{\textsc{BatchLP}} \\
\cmidrule(lr){2-5}\cmidrule(lr){6-9}\cmidrule(lr){10-10}
Instance & $m$ & $n$ & $nnz$ & $|S|$ & 8 threads & 16 threads & 32 threads & 64 threads & GPU B200 \\
\midrule

\multicolumn{10}{l}{\textbf{Combinatorial Auctions}} \\
inst\_100\_1500  & 360  & 1399 & 8418  & 154 & 0.03 & 0.02 & 0.02 & 0.03 & 0.05 \\
inst\_300\_1500  & 557  & 1453 & 8547  & 512 & 0.60 & 0.41 & 0.26 & 0.27 & 0.09 \\
inst\_300\_3000  & 857  & 2840 & 16635 & 514 & 0.73 & 0.48 & 0.31 & 0.36 & 0.08 \\
inst\_500\_1500  & 755  & 1466 & 7989  & 800 & 2.24 & 1.42 & 1.04 & 0.97 & 0.20 \\
inst\_500\_3000  & 1025 & 2867 & 16493 & 846 & 2.51 & 1.51 & 1.24 & 1.04 & 0.09 \\
\multicolumn{5}{r}{\textbf{Total runtime}} & \textbf{6.11} & \textbf{3.84} & \textbf{2.87} & \textbf{2.67} & \textbf{0.51} \\
\multicolumn{5}{r}{\textbf{Speed-up}}      & \textbf{1.0} & \textbf{1.6} & \textbf{2.1} & \textbf{2.3} & \textbf{12.0} \\
\midrule

\multicolumn{10}{l}{\textbf{Set Covering}} \\
inst\_1000r\_1000c\_0.01 & 1000 & 1000 & 10000  & 330 & 1.15 & 0.63 & 0.47 & 0.49 & 0.11 \\
inst\_1000r\_1000c\_0.05 & 1000 & 1000 & 50000  & 136 & 0.44 & 0.27 & 0.19 & 0.20 & 0.06 \\
inst\_1000r\_2000c\_0.01 & 1000 & 2000 & 20000  & 309 & 1.22 & 0.72 & 0.50 & 0.51 & 0.10 \\
inst\_1000r\_2000c\_0.05 & 1000 & 2000 & 100000 & 116 & 0.31 & 0.20 & 0.15 & 0.18 & 0.06 \\
inst\_1000r\_3000c\_0.01 & 1000 & 3000 & 30000  & 265 & 0.77 & 0.54 & 0.37 & 0.37 & 0.09 \\
inst\_1000r\_3000c\_0.05 & 1000 & 3000 & 150000 & 127 & 0.59 & 0.35 & 0.24 & 0.27 & 0.07 \\
inst\_2000r\_1000c\_0.01 & 1978 & 999  & 19781  & 382 & 2.58 & 1.81 & 1.18 & 1.18 & 0.12 \\
inst\_2000r\_1000c\_0.05 & 2000 & 1000 & 100000 & 182 & 1.96 & 1.06 & 0.79 & 0.76 & 0.07 \\
inst\_2000r\_2000c\_0.01 & 2000 & 2000 & 40000  & 449 & 5.96 & 4.16 & 2.69 & 2.57 & 0.15 \\
inst\_2000r\_2000c\_0.05 & 2000 & 2000 & 200000 & 170 & 1.73 & 1.00 & 0.73 & 0.76 & 0.08 \\
inst\_3000r\_1000c\_0.01 & 3000 & 1000 & 30000  & 493 & 8.71 & 5.31 & 3.69 & 3.51 & 0.15 \\
inst\_3000r\_1000c\_0.05 & 3000 & 1000 & 150000 & 164 & 1.67 & 0.98 & 0.77 & 0.72 & 0.07 \\
inst\_3000r\_3000c\_0.01 & 3000 & 3000 & 90000  & 522 & 17.60 & 8.56 & 6.21 & 5.93 & 0.17 \\
inst\_3000r\_3000c\_0.05 & 3000 & 3000 & 450000 & 191 & 4.14 & 2.25 & 1.82 & 1.78 & 0.11 \\
\multicolumn{5}{r}{\textbf{Total runtime}} & \textbf{48.83} & \textbf{27.84} & \textbf{19.80} & \textbf{19.23} & \textbf{1.41} \\
\multicolumn{5}{r}{\textbf{Speed-up}}      & \textbf{1.0} & \textbf{1.8} & \textbf{2.5} & \textbf{2.5} & \textbf{34.6} \\

\midrule

\multicolumn{10}{l}{\textbf{Maximum Independent Set}} \\
inst\_1000\_12 & 8060  & 1000 & 18815 & 1888 & 47.71  & 22.91 & 17.69 & 17.20 & 0.29 \\
inst\_1000\_4  & 8012  & 1000 & 18770 & 1882 & 62.79  & 31.74 & 24.24 & 23.38 & 0.28 \\
inst\_1000\_8  & 8209  & 1000 & 18982 & 1870 & 63.23  & 32.43 & 22.56 & 21.33 & 0.29 \\
inst\_2000\_12 & 18992 & 2000 & 41213 & 3320 & 88.56  & 52.97 & 34.86 & 35.01 & 0.87 \\
inst\_2000\_4  & 18832 & 2000 & 40960 & 3320 & 65.39  & 35.53 & 27.00 & 26.16 & 0.83 \\
inst\_2000\_8  & 18961 & 2000 & 41100 & 3328 & 77.92  & 40.23 & 29.28 & 29.09 & 0.84 \\
inst\_3000\_12 & 35612 & 3000 & 71338 & 5226 & 115.12 & 63.56 & 43.92 & 45.12 & 2.02 \\
inst\_3000\_4  & 35582 & 3000 & 71289 & 5218 & 130.66 & 71.48 & 52.26 & 52.57 & 1.99 \\
inst\_3000\_8  & 35604 & 3000 & 71325 & 5204 & 105.64 & 69.68 & 42.99 & 43.46 & 1.96 \\
\multicolumn{5}{r}{\textbf{Total runtime}} & \textbf{757.02} & \textbf{420.53} & \textbf{294.80} & \textbf{293.32} & \textbf{9.37} \\
\multicolumn{5}{r}{\textbf{Speed-up}}      & \textbf{1.0} & \textbf{1.8} & \textbf{2.6} & \textbf{2.6} & \textbf{80.8} \\
\midrule

\multicolumn{10}{l}{\textbf{Facility Location}} \\
inst\_1200\_200\_10 & 241401 & 240200 & 960400 & 92 & 410.71 & 283.56 & 247.59 & 215.15 & 0.41 \\
inst\_1200\_200\_15 & 241401 & 240200 & 960400 & 88 & 192.56 & 117.60 & 96.49  & 100.83 & 0.40 \\
inst\_1200\_200\_5  & 241401 & 240200 & 960400 & 62 & 135.40 & 110.95 & 88.99  & 81.05  & 0.38 \\
inst\_400\_200\_10  & 80601  & 80200  & 320400 & 40 & 22.52  & 18.09  & 18.05  & 17.97  & 0.12 \\
inst\_400\_200\_15  & 80601  & 80200  & 320400 & 26 & 24.68  & 16.73  & 19.17  & 17.42  & 0.11 \\
inst\_400\_200\_5   & 80601  & 80200  & 320400 & 32 & 4.84   & 4.97   & 5.64   & 5.05   & 0.11 \\
inst\_800\_200\_10  & 161001 & 160200 & 640400 & 58 & 132.33 & 93.22  & 81.69  & 85.90  & 0.24 \\
inst\_800\_200\_15  & 161001 & 160200 & 640400 & 60 & 108.64 & 76.74  & 69.47  & 65.15  & 0.22 \\
inst\_800\_200\_5   & 161001 & 160200 & 640400 & 32 & 33.22  & 23.98  & 30.13  & 24.30  & 0.19 \\
\multicolumn{5}{r}{\textbf{Total runtime}} & \textbf{1064.90} & \textbf{745.84} & \textbf{657.22} & \textbf{612.82} & \textbf{2.18} \\
\multicolumn{5}{r}{\textbf{Speed-up}}      & \textbf{1.0} & \textbf{1.4} & \textbf{1.6} & \textbf{1.7} & \textbf{488.5} \\
\bottomrule
\end{tabular}
\end{table*}

\paragraph{Implementation details.}
We implemented the code in standard C++17 and used the CUDA Toolkit and the Thrust library to develop the main GPU kernels.
Our implementation does not perform any dynamic memory allocation at runtime. All necessary data is transferred to the GPU before use according to the techniques described in Section~\ref{sec:batchLPs}.

For the matrix operations, we rely on the {cuSPARSE} and {cuBlas} libraries. We store in memory both spare matrices $A$ and $A^\tpose$, and we use {\tt cusparseSpMM} to compute $A X$ and $A^\tpose Y$, that is, the product between sparse and dense matrices.We use {\tt cublasDgemv} to compute matrix vector operations such as, for instance, $X^\tpose c$, which compute $c^\tpose x_i$ for every column $x_i$ of $X$, and {\tt cublasDnrm2} for vector norms. All other functions are executed on the GPU using our own customized CUDA kernels.

\paragraph{Hardware details.}
We use a B200 and a H100 on the {\tt brev} NVIDIA cloud infrastructure mounted on a virtual machine equipped with an Intel Xeon CPU with 56 physical cores and 128 GB of RAM.

\begin{figure}[t]
\centering
\begin{tikzpicture}
  \begin{axis}[
    ybar,
    bar width=6pt,
    width=\ifpreprint 0.7\textwidth\else\columnwidth\fi,
    height=5.5cm,
    ymode=log,
    log origin=infty,
    symbolic x coords={CA, SC, MIS, FL},
    xtick=data,
    xticklabels={Comb.\ Auct., Set Cover., Max Ind.\ Set, Facility Loc.},
    x tick label style={font=\small},
    ylabel={Speedup vs.\ 8 threads},
    ylabel style={font=\footnotesize},
    y tick label style={font=\footnotesize},
    ymin=1,
    ymax=1000,
    enlarge x limits=0.15,
    legend style={
      at={(0.02,0.98)},
      anchor=north west,
      legend columns=1,
      font=\small,
      draw=none,
      fill=none,
      cells={anchor=west},
    },
    area legend,
  ]
  \addplot[fill=black!20, draw=black!40] coordinates {(CA,1.6) (SC,1.8) (MIS,1.8) (FL,1.4)};
  \addplot[fill=black!40, draw=black!60] coordinates {(CA,2.1) (SC,2.5) (MIS,2.6) (FL,1.6)};
  \addplot[fill=black!60, draw=black!70] coordinates {(CA,2.3) (SC,2.5) (MIS,2.6) (FL,1.7)};
  \addplot[fill=black, draw=black] coordinates {(CA,12.0) (SC,34.6) (MIS,80.8) (FL,488.5)};
  \legend{16 threads, 32 threads, 64 threads, GPU}
  \end{axis}
\end{tikzpicture}
\caption{Speedup of FSB runtime vs.\ 8-thread CPU baseline. GPU achieves 12--489$\times$ speedup. Full results in Table~\ref{tab:cuopt:threads}.}
\label{fig:fsb-summary}
\end{figure}
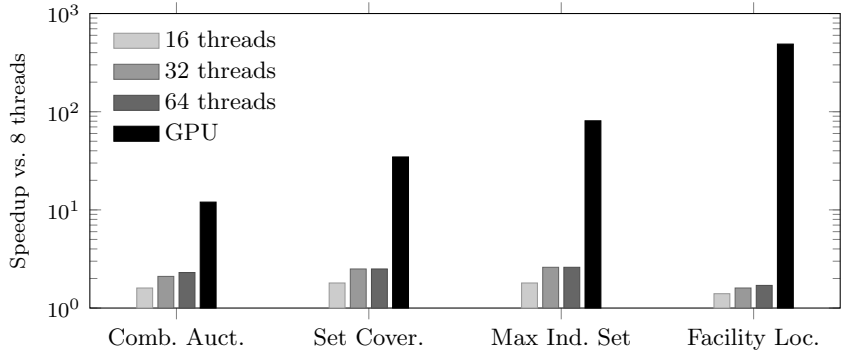
\subsection{FSB integration within {\sc cuOpt} MIP solver}
We integrate \textsc{BatchLP} into the branch-and-bound method used in NVIDIA {\sc cuOpt} 26.02 \cite{cuopt}.
{\sc cuOpt} performs a round of strong branching after solving the root relaxation to initialize pseudocosts for fractional variables.
These pseudocosts are then used to make branching decisions.
{\sc cuOpt}'s default strong-branching implementation splits the subproblems evenly across a set of threads.
Each thread uses dual simplex, warm-started from the root's optimal basis, to solve each LP to optimality.

In our main experiments, we use the instances introduced in \cite{Gasse2019} to perform a computational study on learning branching strategies that can be competitive with full strong branching, while being computable in a very short time.
Here, we take a different perspective: rather than approximating FSB, we aim to execute it by exploiting the GPU, solving each subproblem LP to optimality as fast as possible with {\sc BatchLP}.

Figure~\ref{fig:fsb-summary} summarizes the results; full details are reported in Table~\ref{tab:cuopt:threads}.
Table~\ref{tab:cuopt:threads} compares the time required to perform a single round of strong branching using dual simplex on multiple threads (8, 16, 32, and 64 threads on a Dell workstation with 56 physical cores -- Intel Xeon CPU) versus \textsc{BatchLP}.
We use the random generator from \cite{Gasse2019} to generate instances from four families of combinatorial optimization problems: Combinatorial Auctions, Set Covering, Maximum Independent Set, and Facility Location.
Within each family, instances are ordered by increasing size in terms of the number of constraints $m$, variables $n$, and nonzeros $nnz$.
The number of subproblems to be solved, denoted by $|S|$, varies across problem classes.

For each family, we report the total speedup relative to the runtime of dual simplex in parallel on 8 threads.
The speedup from using more threads is not linear, and the difference between 32 and 64 threads is minimal.
In contrast, {\sc BatchLP} yields a substantial speedup across all instances (except the smallest one); for the Facility Location MIPs, it achieves a speedup of $488.5\times$.


\subsection{OBBT results}
In this section, we use \textsc{BatchLP} independently of cuOpt and customized to perform OBBT by carefully handling the incumbent best dual bound on the objective value. We compare its performance to an OBBT implementation based on the commercial solver Gurobi.
While in strong branching, the lower numerical accuracy of first-order methods is less critical (branching scores are used only as a heuristic to order variables), in OBBT, we must guarantee the use of a \emph{safe} dual bound before restricting the domain of any variable.
For this reason, when running OBBT in \textsc{BatchLP}, we enforce a higher numerical accuracy when checking dual feasibility ($\epsilon_{\text{dual}}=10^{-8}$), and we subtract/add
\[
\Delta = \epsilon \Bigl(1 + |c^\tpose x| + \bigl|\suppf_{[\underline{x}, \bar{x}]}(r) + \suppf_{[l, u]}(y)\bigr|\Bigr),
\]
to the lower/upper bounds computed via OBBT.
Finally, we update a bound only if the improvement is larger than $10^{-4}$.

\begin{figure}[t]
\centering
\begin{tikzpicture}
  \begin{axis}[
    ybar,
    bar width=6pt,
    width=\ifpreprint 0.85\textwidth\else 1.05\columnwidth\fi,
    height=5.5cm,
    ymode=log,
    symbolic x coords={120,232,233,428,438,511,523,53,561,584,604,633,64,698,704,728,744,751,772,791},
    xtick=data,
    x tick label style={font=\scriptsize, rotate=45, anchor=east},
    xlabel={Instance (test\_*)},
    xlabel style={font=\footnotesize},
    ylabel={Speedup vs.\ Gurobi dual simplex},
    ylabel style={font=\footnotesize},
    y tick label style={font=\footnotesize},
    ymin=0.1,
    ymax=2000,
    extra y ticks={1},
    extra y tick style={grid=major, grid style={black, thick}},
    enlarge x limits=0.03,
    area legend,
  ]
  \addplot[fill=black, draw=black] coordinates {
    (120,25.3) (232,16.4) (233,19.0) (428,25.6) (438,21.9)
    (511,18.7) (523,12.1) (53,28.5) (561,43.0) (584,15.2)
    (604,25.3) (633,34.9) (64,25.6) (698,26.9) (704,32.6)
    (728,28.9) (744,26.7) (751,43.6) (772,11.2) (791,25.3)
  };
  \end{axis}
\end{tikzpicture}
\caption{OBBT speedup of \textsc{BatchLP} vs.\ Gurobi dual simplex on neural network verification instances. Average speedup: 25.7$\times$. Full results in Table~\ref{tab:obbt}.}
\label{fig:obbt-summary}
\end{figure}
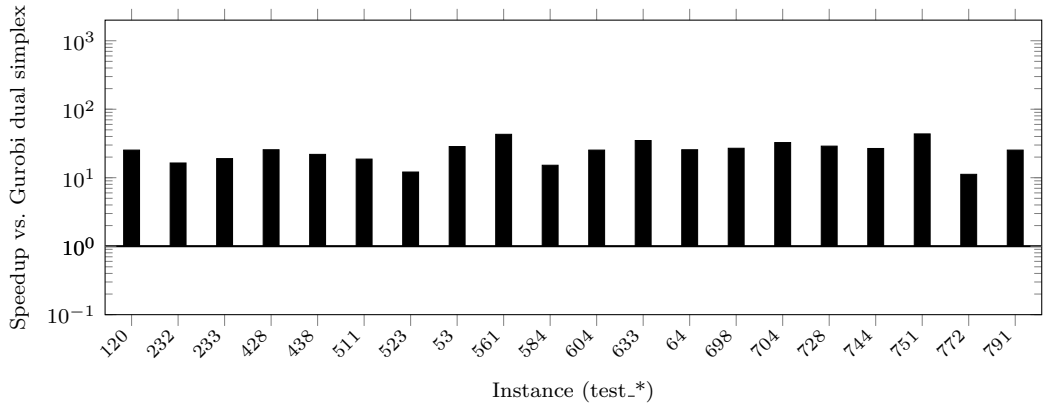

In the numerical experiments, we apply OBBT to neural network verification instances from~\cite{Nair2020}.
Figure~\ref{fig:obbt-summary} summarizes the results, with full details in Table~\ref{tab:obbt}.
Table~\ref{tab:obbt} reports our results, comparing \textsc{BatchLP} with a sequential version of Gurobi. Similarly to Table 1: $m$ is the number of constraints, $n$ is the number of variables, while $|S| = 2n$ is not explicitly reported. Finally, $den$ is the density of the coefficient matrix $A$ (similarly to the $nnz$ value reported in Table \ref{tab:cuopt:threads}).
For Gurobi dual simplex, we report the total runtime and the average number of pivots required to solve each subproblem to optimality.
For \textsc{BatchLP}, we report the total runtime, the speed-up over the sequential version, and the number of subproblems solved before hitting the iteration limit.
The last two columns report the number of variables whose domain changes, and the average domain reduction (as a percentage of the original domain).
We remark that, since we perform OBBT on the model obtained after Gurobi MIP presolve, several standard bound-reduction techniques have already been applied to the original model.

\paragraph{Interpreting the speed-up.}
The results in Table \ref{tab:obbt} show a large sped-up (25.7 on average) that is significant considering the quality of the Gurobi simplex implementation. Such an average speed-up of \textsc{BatchLP} with respect to sequential dual simplex should be interpreted as follows: under ideal CPU parallelism (which we have shown in Table \ref{tab:cuopt:threads} does not happen), solving in parallel the strong branching LPs with dual simplex using 26 threads on 26 physical cores should yield a runtime comparable to that of {\sc BatchLP} running on a single GPU.

\subsection{Results on MIPLIB 2017}
Finally, we present preliminary results on running {\sc BatchLP} on MIPLIB 2017 using the default strong-branching strategy we implemented in cuOpt 26.02, i.e., running the dual simplex in parallel, with a limit of only 200 pivots per subproblem.
Figure~\ref{fig:miplib-summary} summarizes the results; full details are reported in Table~\ref{tab:miplib} in the Appendix.
Table~\ref{tab:miplib} reports results for a subset of instances on which {\sc BatchLP} on an H100 GPU is competitive with the dual simplex running in parallel on 28 threads (i.e., 28 physical cores).
Unsurprisingly, for all other MIPLIB 2017 instances not reported in the table, the dual simplex on 28 threads (with the 200-pivot limit) is faster than {\sc BatchLP}.

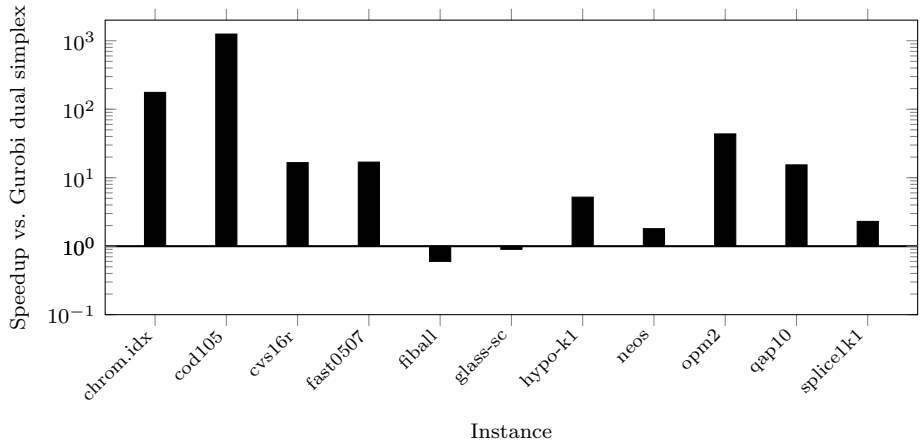
\begin{figure}[t]
\centering
\begin{tikzpicture}
  \begin{axis}[
    ybar,
    bar width=8pt,
    width=\ifpreprint 0.75\textwidth\else 1.05\columnwidth\fi,
    height=5.5cm,
    ymode=log,
    symbolic x coords={chrom,cod105,cvs16r,fast,fiball,glass,hypo,neos,opm2,qap10,splice},
    xtick=data,
    xticklabels={chrom.idx, cod105, cvs16r, fast0507, fiball, glass-sc, hypo-k1, neos, opm2, qap10, splice1k1},
    x tick label style={font=\scriptsize, rotate=45, anchor=east},
    xlabel={Instance},
    xlabel style={font=\footnotesize},
    ylabel={Speedup vs.\ Gurobi dual simplex},
    ylabel style={font=\footnotesize},
    y tick label style={font=\footnotesize},
    ymin=0.1,
    ymax=2000,
    extra y ticks={1},
    extra y tick style={grid=major, grid style={black, thick}},
    enlarge x limits=0.07,
    area legend,
  ]
  \addplot[fill=black, draw=black] coordinates {
    (chrom,176.2) (cod105,1251.7) (cvs16r,16.6) (fast,16.9) (fiball,0.6)
    (glass,0.9) (hypo,5.2) (neos,1.8) (opm2,43.6) (qap10,15.4) (splice,2.3)
  };
  \end{axis}
\end{tikzpicture}
\caption{MIPLIB 2017 speedup of \textsc{BatchLP} vs.\ cuOpt dual simplex (28 threads). Values below 1 indicate dual simplex is faster. Full results in Table~\ref{tab:miplib}.}
\label{fig:miplib-summary}
\end{figure}

Although the impact of moving strong branching to the GPU has been presented only on few MIP instances and within a solver that is admittedly less developed than the commercial ones like Gurobi, it is important to note that the integration of such GPU mechanism within a state-of-the-art solver is a standalone research question that deserves specific work beyond the scope of the current paper. Initial tests with cuOpt 26.02 on a family of `qap' instances (see instance qap10 in Table \ref{tab:miplib}) where we run the entire branch and bound using either {\sc BatchLP} or the default parallel simplex implementation show very promising results but they are too preliminary to be reported in full.

\section{Conclusions}
\label{sec:conclusions}
In this paper, we introduced a batched first-order method for solving multiple linear programs in parallel on GPUs, specifically focusing on mixed-integer programming (MIP) challenges like strong branching and optimization-based bound tightening. Our approach extends the primal-dual hybrid gradient algorithm, leveraging matrix-matrix operations to fully exploit the capabilities of modern GPU architectures. Through a series of computational experiments, we demonstrated significant speed-ups over traditional simplex-based methods, particularly in instances with large problem sizes and numerous subproblems in the batches. 
This advancement offers the MIP community a powerful tool for handling repeated subproblem solutions more efficiently, paving the way for fully incorporating GPUs in MIP algorithms. 




\ifpreprint
\section*{Acknowledgments}
\myack
\else\ifdefined\isaccepted
\section*{Acknowledgments}
\myack
\fi\fi

\ifpreprint\else
\section*{Impact Statement}
This paper presents work whose goal is to advance the field of mathematical optimization.
Our methods improve the computational efficiency of solving mixed-integer programs by exploiting GPU parallelism, which may reduce energy consumption for large-scale optimization tasks.
There are many potential societal consequences of our work, none of which we feel must be specifically highlighted here.
\fi

\ifpreprint
\bibliography{bibliography}
\else
\balance
\bibliography{bibliography}
\bibliographystyle{icml2026}
\fi

\newpage
\appendix
\onecolumn


\section{OBBT Results}
\label{app:obbt-results}
Table~\ref{tab:obbt} presents detailed results for the OBBT experiments discussed in Section~\ref{sec:results}.

\begin{table*}[h]
\centering
\caption{Results for OBBT applied after MIP presolve of Gurobi: Dual Simplex--Gurobi 13.0 on Intel Xeon Platinum 8592; {\sc BatchLP} iteration limit $100\,000$ (always hit) on a GPU B200 (cost of machine 5.29 USD/hour).}
\label{tab:obbt}
\small
\setlength{\tabcolsep}{4pt}
\begin{tabular}{lrrr|rr|rrc|rr}
\toprule
& \multicolumn{3}{c}{} & \multicolumn{2}{|c}{Dual simplex} & \multicolumn{3}{|c}{{\sc BatchLP}} & \multicolumn{2}{|c}{var changed} \\
Instance & $m$ & $n$ & den & runtime & pivots & runtime & Speed-up & solved & num & perc. \\
\midrule
test\_120 & 1564 & 1767 & 1.8\% & 507.9 & 755.9 & 20.1 & 25.3 & 3429/3534 & 141 & 4.0\% \\
test\_232 & 1296 & 1491 & 2.0\% & 271.2 & 594.1 & 16.6 & 16.4 & 2866/2982 & 105 & 3.5\% \\
test\_233 & 1315 & 1512 & 1.9\% & 225.8 & 495.6 & 11.9 & 19.0 & 2973/3024 & 110 & 3.6\% \\
test\_428 & 1485 & 1617 & 2.0\% & 551.7 & 887.6 & 21.6 & 25.6 & 2859/3234 & 93 & 2.9\% \\
test\_438 & 1229 & 1332 & 2.3\% & 263.7 & 660.7 & 12.0 & 21.9 & 2591/2664 & 112 & 4.2\% \\
test\_511 & 1100 & 1208 & 2.8\% & 250.8 & 734.4 & 13.4 & 18.7 & 1800/2416 & 168 & 7.0\% \\
test\_523 & 932 & 1036 & 3.1\% & 145.0 & 622.4 & 11.9 & 12.1 & 1981/2072 & 134 & 6.5\% \\
test\_53  & 1635 & 1741 & 1.6\% & 485.3 & 760.3 & 17.1 & 28.5 & 3287/3482 & 108 & 3.1\% \\
test\_561 & 1457 & 1643 & 2.0\% & 573.8 & 902.9 & 13.4 & 43.0 & 3167/3286 & 134 & 4.1\% \\
test\_584 & 883 & 980 & 3.2\% & 141.2 & 666.9 & 9.3 & 15.2 & 1922/1960 & 117 & 6.0\% \\
test\_604 & 1362 & 1529 & 2.1\% & 385.5 & 719.3 & 15.3 & 25.3 & 2730/3058 & 59 & 1.9\% \\
test\_633 & 1453 & 1614 & 1.9\% & 462.0 & 786.2 & 13.2 & 34.9 & 3119/3228 & 150 & 4.7\% \\
test\_64  & 1589 & 1736 & 1.7\% & 503.6 & 766.2 & 19.7 & 25.6 & 3290/3472 & 108 & 3.1\% \\
test\_698 & 1347 & 1494 & 2.0\% & 309.5 & 665.9 & 11.5 & 26.9 & 2915/2988 & 115 & 3.9\% \\
test\_704 & 1574 & 1731 & 1.7\% & 513.5 & 769.3 & 15.7 & 32.6 & 3369/3462 & 129 & 3.7\% \\
test\_728 & 1410 & 1570 & 2.0\% & 393.1 & 710.2 & 13.6 & 28.9 & 3004/3140 & 121 & 3.9\% \\
test\_744 & 1575 & 1801 & 1.6\% & 435.6 & 662.8 & 16.3 & 26.7 & 3513/3602 & 124 & 3.4\% \\
test\_751 & 1672 & 1927 & 1.6\% & 613.6 & 784.3 & 14.1 & 43.6 & 3751/3854 & 132 & 3.4\% \\
test\_772 & 907 & 973 & 3.2\% & 128.9 & 633.0 & 11.5 & 11.2 & 1572/1946 & 99 & 5.1\% \\
test\_791 & 1410 & 1548 & 2.1\% & 432.0 & 773.0 & 17.1 & 25.3 & 2899/3096 & 100 & 3.2\% \\
\midrule
mean &  &  &  & 379.7 &  & 14.8 & 25.7 &  &  & 4.1\% \\
\bottomrule
\end{tabular}
\end{table*}

\section{MIPLIB 2017 Results}
\label{app:miplib-results}
Table~\ref{tab:miplib} presents detailed results for the MIPLIB 2017 experiments discussed in Section~\ref{sec:results}.

\begin{table*}[h]
\centering
\caption{Comparison of cuOpt dual simplex running on 28 threads (28 physical cores) with {\sc BatchLP} running on a GPU H100.\label{tab:miplib}}
\small
\setlength{\tabcolsep}{5pt}
\begin{tabular}{lrrrrrrr}
\toprule
Instance & $m$ & $n$ & $nnz$ & $|S|$ & Dual simplex  & \textsc{BatchLP} & Speed-up \\
\midrule
chromaticindex512-7      & 33791  & 36864 & 135156  & 51604 & 3522.4 & 20.0   & 176.2 \\
cod105                   & 1024   & 1024  & 57344   & 256   & 212.8  & 0.2    & 1251.7 \\
cvs16r128-89             & 4633   & 3472  & 12528   & 6420  & 20.9   & 1.3    & 16.6 \\
fast0507                 & 482    & 62171 & 401756  & 556   & 5.8    & 0.3    & 16.9 \\
fiball                   & 2387   & 32899 & 101452  & 518   & 22.8   & 37.1   & 0.6 \\
glass-sc                 & 6119   & 214   & 63918   & 202   & 1.0    & 1.1    & 0.9 \\
hypothyroid-k1           & 5189   & 2595  & 431326  & 4806  & 3292.7 & 627.9  & 5.2 \\
neos-5052403-cygnet      & 19134  & 27593 & 2448853 & 1272  & 63.7   & 35.8   & 1.8 \\
opm2-z10-s4              & 146325 & 5958  & 340288  & 11170 & 1921.5 & 44.1   & 43.6 \\
qap10                    & 1820   & 4150  & 18200   & 2442  & 28.8   & 1.9    & 15.4 \\
splice1k1                & 6504   & 3252  & 1758012 & 3816  & 3355.3 & 1468.6 & 2.3 \\
\bottomrule
\end{tabular}
\end{table*}

\end{document}